\documentclass[12pt]{amsart}

\usepackage{amssymb,amsmath,latexsym,amsthm}
\usepackage[english]{babel}

\newtheorem{theorem}{Theorem}[section]

\theoremstyle{definition}

\begin{document}

\title{Selectivity on division algebras}


\author{\sc Luis Arenas-Carmona}


\newcommand\Q{\mathbb Q}
\newcommand\reali{\mathbb R}
\newcommand\compleji{\mathbb C}
\newcommand\alge{\mathfrak{A}}
\newcommand\Da{\mathfrak{D}}
\newcommand\Ha{\mathfrak{H}}
\newcommand\oink{\mathcal O}
\newcommand\matrici{\mathbb{M}}
\newcommand\Txi{\lceil}
\newcommand\ad{\mathbb{A}}
\newcommand\enteri{\mathbb Z}
\newcommand\finitum{\mathbb{F}}
\newcommand\bbmatrix[4]{\left(\begin{array}{cc}#1&#2\\#3&#4\end{array}\right)}

\maketitle

\begin{abstract}
A commutative order in a central simple algebra over a number field is said to be selective if it embeds in some,
 but not
all, maximal orders in the algebra. We completely characterize selective orders in central division algebras, 
of dimension $9$ or greater, in terms of the characterization of selective orders given by Chindburg and Friedman 
in the quaternionic case.
\end{abstract}

\bigskip
\section{Introduction}

Back in 1936, Chevalley \cite{Chevalley} proved that the number $q$ of conjugacy classes of maximal orders, in a matrix algebra $B$, over a number field $K$, 
containing a copy of the ring of integers of a fixed maximal  subfield $L$ of $B$ had the form
$$q=N/[H\cap L:K],$$ where $N$ is the total number of conjugacy classes and
 $H$ is the Hilbert class field of $K$. An analog result was later proved in 1999 by Chinburg
and Friedman \cite{FriedmannQ} for a quaternion algebra $Q$ satisfying the Eichler condition, where they gave 
necessary and sufficient conditions for  a commutate suborder to be contained in every maximal order in $Q$. 
When this fails to be the case, the order is called selective and embeds into exactly one half of all conjugacy
classes. This result was later extended to Eichler orders by several authors \cite{Guo}, \cite{Chan}.
We ourselves proved in \cite{spinor} that Chevalley's formula $q=N/[H\cap L:K]$ remains valid for central simple
algebras without partial ramification, i.e., for algebras that are either division algebras,
or pure matrix algebras, at every finite place, as long as the Hilbert class field is replaced by another class field,
which we called the spinor class field $\Sigma$. Furthermore  in 2012, in a second paper  \cite{abelianos}, 
we proved we always have the identity $q=N/[F:K]$, for some explicit subfield $F\subseteq L\cap\Sigma$ defined
by certain splitting conditions, although the equality $F=L\cap\Sigma$ does not hold in general. The field $F$
is called the representation field.

In all of this work, a commutative order is called selective for a central simple algebra $B$ if it embeds in some, but not
all, the conjugacy classes of maximal orders in $B$. In the notations of the previous paragraph, an order $\Ha$
is selective if and only if the representation field $F=F(\Ha)$ strictly contains $K$. 
The same year that \cite{abelianos} was published, Linowitz and Shemanske \cite{lino2} noted that selectivity
was not possible in a division algebra $B$ when $\sqrt{\mathrm{dim}_K(B)}$ is an odd prime. Furthermore,
in a recent pre-print \cite{lino3}, the same authors proved that selectivity is not possible in an algebra 
$B$ that had a fully ramified 
place. The condition implies, in particular, that $B$ is a division algebra. This hinted us that the situation for division 
algebras could be simpler than for more general central simple algebras. While it is easy to find commutative orders
of rank $n$, in a $n^2$ dimensional matrix algebra, whose the selectivity rate $q$ is $1/n$, provided that $n$ divides the annihilator of the abelian group $\mathrm{Gal}(H/K)$, the corresponding rate in division algebras can be only $0$, $1$, or
$1/2$.  In fact, we show in Theorem 1 below that the selectivity
that was found by Chindburg and Friedman in quaternion algebras is, in a sense, responsible for all selectivity on division algebras.

\begin{theorem}
Let $K$ be a number field. Let $B$ be a division $K$-algebra, and let $L\subseteq B$ be a maximal subfield. Let $\Ha$ be
an arbitrary order in $L$. If $\Ha$ is selective in $B$, the following conditions hold:
\begin{enumerate}
\item $\sqrt{\mathrm{dim}_K(B)}$ is twice an odd number,
\item the ramification degree of $B$ at every finite place is odd,
\item $B$ has a decomposition $B=B_1\otimes_K B_2$, where $B_1$ is a quaternion algebra and $B_2$ is an odd dimensional division algebra. 
\item $L$ contains a quadratic extension that is selective for $B_1$.
\end{enumerate} 
On the other hand, if $L$ and $B$ satisfy these conditions, then $\oink_L$ embeds into exactly one
 half of all isomorphism classes of maximal  orders in $B$.
\end{theorem}

Is it worth to note that, aside from the last statement, the proof of this result uses only the trivial 
bound $F\subseteq L\cap\Sigma$, and not the
precise computation of $F$ in \cite{abelianos}. The main idea behind the proof, which also seems to be present in the
previous work by Linowitz and Shemanske, is that a field that intersects "too much" the spinor class field $\Sigma$ must forfeit the necessary conditions to be embedded in the algebra $B$.

    The arithmetic of maximal orders in quaternion division algebras is central in the construction of 
isospectral non-isometrical  varieties in \cite{vigneras2} and \cite{Rajan}.
 Lack of selectivity is particularly important in this construction see for example \cite[\S5.1]{linothesis}.
 
\section{The brauer group of a number field}
In this section we recall the structure results for the Brauer group of a number field that are needed to prove our
 main result. The proof of these results can be found in \cite{albert}. Recall that every central simple algebra over a field $K$ is isomorphic to a ring of matrices $\matrici_n(D)$ with entries in a central division $K$-algebra $D$. 
The Brauer group $\mathrm{Br}(K)$
of $K$ is a group whose element are equivalence classes of central simple algebras over $K$, where two algebras
$\matrici_n(D)$ and  $\matrici_m(B)$ are in the same class if and only if the division algebras $D$ and $B$ are isomorphic.
In this group, the identity is the class of the $K$-algebra $K$, and the inverse of a class[C] is the class $[C^{\mathrm{op}}]$ corresponding to the algebra $C^{\mathrm{op}}$ which is isomorphic to $C$ as a vector space,
but endowed with the opposite multiplication $a*b=b\cdot a$. The group operation is $[C][D]=[C\otimes_K D]$.

For any field $F$ containing $K$, there is a homomorphism $\phi:\mathrm{Br}(K)\rightarrow\mathrm{Br}(F)$ sending
the class of a $K$-algebra $C$ to the class of the central simple $F$-algebra $C\otimes_KF$. In these terms, for any number field $K$ there exists an exact sequence
$$0\rightarrow \mathrm{Br}(K)\rightarrow\prod_{\wp\in\Pi(K)}\mathrm{Br}(K_\wp)\rightarrow\mathbb{Q}/\enteri
\rightarrow0,$$
where $\Pi(K)$ is the set of all places of $K$, including the archimedean places.
In particular, a central simple $K$-algebra $\alge$ is completely determined by the set of completions $\alge_\wp$.
The brauer group of the completion $K_\wp$ is given by
$$\mathrm{Br}(K_\wp)\cong\left\{\begin{array}{ll}\enteri/2\enteri&\textnormal{if }K_\wp\cong\reali\\
\{0\}&\textnormal{if }K_\wp\cong\compleji\\ \mathbb{Q}/\enteri&\textnormal{otherwise}\end{array}\right..$$
Identifying $\enteri/2\enteri$ with $\left(\frac{1}{2}\enteri\right)/\enteri$, one can associate to a central simple algebra $\alge$, the Hasse invariant $i_\wp(\alge)\in\mathbb{Q}/\enteri$, for every local place $\wp$,
 with the obvious restrictions at the archimedean places, and satisfying the properties
$$i_\wp(\alge)=0\textnormal{ for almost all }\wp,\quad \sum_{\wp\in\Pi(K)}i_\wp(\alge)=0.$$
Furthermore, for any positive integer $n$, the elements of order $n$ in the Brauer
group are exactly the clases of $n^2$-dimensional central division $K$-algebras. It follows that every element in $\mathrm{Br}(K)$
annihilated by $n$ is the class of a $n^2$-dimensional central simple $K$-algebra. The order of any class $[C]$
in the brauer group $\mathrm{Br}(K)$ is the least common multiple of the orders (i.e., denominators) of the
invariants $i_\wp(C)$. This local order is, by definition, the ramification degree $e_\wp(B/K)$.
In particular, if $C$ is an $n^2$-dimensional division algebra, and $p$ is a prime
dividing $n$, there must exists a local place $\wp$ such that $p^{v_p(n)}$ divides $e_\wp(B/K)$,
where $v_p(n)$ is the $p$-adic valuation of $n$, i.e., the largest integer $v$ such that $p^v$ divides $n$.
 This observation is used in all that follows.

\section{The theory of representation fields}

Let $K$ be a number field and let $\alge$ be a central simple $K$-algebra. Two orders of maximal rank $\Da$ and $\Da'$
are said to be in the same genus if the completions $\Da_\wp$ and $\Da'_\wp$ are conjugate for every finite place $\wp\in\Pi(K)$. For simplicity, we write $\Da_\wp=\alge_\wp$ whenever $\wp\in\Pi(K)$ is archimedean.  If we define the adelizations $\alge_\ad$ and $\Da_\ad$ of the algebra $\alge$ and the
order $\Da$ by the formulas $\Da_\ad=\prod_{\wp\in\Pi(K)}\Da_\wp$ and\footnote{
It is well known that the definition of $\alge_\ad$ given here is independent of the choice of $\Da$, see fon example \cite{abelianos}.}
$$\alge_\ad=\left\{a\in\prod_{\wp\in\Pi(K)}\alge\wp\Big|a_\wp\in\Da_\wp\textnormal{ for almost every }\wp\in\Pi(K)\right\},$$
two orders $\Da$ and $\Da'$ are in the same genus if and only if $\Da'_\ad=a\Da_\ad a^{-1}$ for some element
$a\in\alge_\ad^*$. They are in the same spinor genus when $a$ can be chosen in the form $a=bc$ where
every coordinate $b_\wp$ of $b$ has reduced norm $1$, while $c\in\alge$ is a global element. It is known that
two orders in the same spinor genus are conjugate whenever the Eichler Condition holds, i.e., $\alge$ is not the algebra
of Hamilton's Quaternions for at least one infinite place $\wp$. This is certainly the case when the dimension of $\alge$
is not $4$, or when $K$ is not totally real. On the other hand, the set of spinor genera in a genus
$\mathrm{gen}(\Da)$ is in correspondence with the Galois group over $K$ 
of the spinor class field $\Sigma=\Sigma(\Da)$, i.e., the class field corresponding to the class group $K^*H(\Da)$
where \begin{equation}\label{uno}
H(\Da)=\{N(a)|a\in\alge_\ad^*, a\Da_\ad a^{-1}=\Da_\ad\}\subseteq J_K.\end{equation} 
The spinor genus containing a second order $\Da'$ corresponds to the element 
$\rho(\Da,\Da')\in\mathrm{Gal}(\Sigma/K)$ defined by $\rho(\Da,\Da')=[N(a),\Sigma/K]$ where $\Da'_\ad=a\Da_\ad a^{-1}$, and $x\mapsto[x,\Sigma/K]$ is the Artin map on ideles. In particular, the number of spinor genera containing
some supra-order $\Da'\in\mathrm{gen}(\Da)$ of a fixed order $\Ha\subseteq\alge$ is $[\Sigma:F]$, where $F$ is
the class field corresponding to the class group $K^*H(\Da|\Ha)$,
where $$H(\Da|\Ha)=\{N(a)|a\in\alge_\ad^*, a\Ha_\ad a^{-1}\subseteq\Da\},$$ provided that  $K^*H(\Da|\Ha)$
is actually a group. This is not always the case, but it holds when $\Ha$ is commutative and $\Da$ is maximal
\cite{abelianos}. It is immediate
from the definition of $H(\Da|\Ha)$, that we have the contention $F\subseteq L\cap\Sigma$ for any maximal field $L$ containing $\Ha$. This is the only bound needed in the sequel. A stronger bound of the form  $F\subseteq F_0(\alge|\Ha)\cap\Sigma$, was used in \cite{continuity} to prove that an order contained in a quadratic extension is spinor selective for some genus in exactly one quaternion algebra. The field $F_0(\alge|\Ha)$ is the class field corresponding to a class group $K^*H(\Ha)$ whose
definition is analogous to (\ref{uno}).  

There is a unique genus of maximal orders in an $n^2$-dimensional central simple algebra $\alge$. Furthermore, the spinor class field $\Sigma_0$ for the genus of maximal orders is characterized as
the largest unramified exponent-$n$ extension of the ground field $K$ satisfying the following conditions \cite[\S2]{spinor}:
\begin{itemize}
 \item $\Sigma_0/K$ is unramified at every infinite place $\wp$ where $\alge_\wp$ is a matrix algebra.
\item At every finite place $\wp$, the inertia degree $f_\wp(\Sigma_0/K)$ divides $f=f_\wp(\alge/K)$, where $\alge_\wp\cong\matrici_f(\mathfrak{B})$ for a local division algebra $\mathfrak{B}$.
\end{itemize}

\section{Proof of the main result}

Let $B$ be an $n^2$-dimensional central division $K$-algebra, and let $p$ be a prime dividing $n$. 
In particular,$B$ has order $n$ in the Brauer group of $K$ (\S2). 
Furthermore, for every fixed prime $p$, there exists a 
finite place $\wp$ such that the ramification degree $e_\wp(B/K)$ is divisible by $p^t$, where $t=v_p(n)$ 
is the $p$-adic valuation of $n$.  By the description of the spinor class field given in \S3, the inertia degree 
$f_\wp(\Sigma/K)$ must be relatively prime to $p$, and the same holds for $f_\wp(F/K)$, where $F$ is the 
representation field. Recall also that the ramification degree $e_\wp(F/K)$ is $1$ since $\Sigma/K$ is 
unramified at finite places.

Now let $\mathbb{P}$ be a place of $L$ lying over $\wp$. The condition that $L$ embeds into $B$ implies that $e_\wp(B/K)$  divides the local degree $[L_{\mathbb{P}}:K_\wp]$. In particular, $p^t$ divides $[L_{\mathbb{P}}:K_\wp]$. Since the extension $F/K$ is abelian, we conclude that
$p^t$ divides $[L_{\mathbb{P}}:F_{\wp'}]$ where $\wp'$ is the restriction of $\mathbb{P}$ to $F$.
 Since this holds for every place $\mathbb{P}$ of $L$ lying over $\wp$, we conclude
that  $p^t$ divides
$$[L:F]=\sum_{\mathbb{P}|\wp'}[L_{\mathbb{P}}:F_{\wp'}].$$
In particular, $p$ does not divide $[F:K]$. 

A similar argument holds for $p=2$ when $v_2(n)\geq2$. This proves (1). On the other hand, if $v_2(n)=1$ and  the ramification degree $e(B_\wp/K_\wp)$ is even, for a finite place $\wp$, the same argument proves that $2$ does not 
divide $[F:K]$ and therefore again $[F:K]=1$, so (2) follows. 

We assume (1) and (2) in all that follows. In particular, $B$ has odd ramification degree at every finite prime, 
while the set $T$, of real places where $B$ has ramification degree $2$, is non-empty. The order of the set $T$ 
must be even, since $\sum_\wp i_\wp(B)\in\mathbb{Q}/\mathbb{Z}$, and therefore also its projection on the $2$-torsion 
of $\mathbb{Q}/\mathbb{Z}$, must vanish. There exists, therefore, a quaternion algebra $B_1$ ramifying at exactly 
the places in $T$. Similarly, there exists an odd dimensional CSA $B_2$ ramifying at the same finite places as $B$ 
and with the same Hasse invariants, satisfying $\mathrm{dim}_K(B_2)=n/2^{v_2(n)}$. Note
that $B_2$ is a division algebra since, for every fixed odd prime $p$ dividing $n$, there exists a finite place 
$\wp$ such that the ramification degree $e_\wp(B_1/K)$ is divisible by $p^{v_p(n)}$. 
Furthermore, both $B$ and $B_1\otimes_K B_2$ have the same invariants. We conclude (3).

Recall that the proportion of spinor genera representing $\Ha$, in the genus $\mathbb{O}_0$ of maximal orders,
is $[F:K]^{-1}$. The first part of the argument showed that $[F:K]=2$, since we are assuming the order is selective 
and every prime-power divisor larger than $2$ has been ruled out as a divisor of $[F:K]$. 
In particular $F$ is a quadratic extension of $K$. 
We claim that it satisfies the conditions in (4).

First we prove that $F\subseteq L$ embeds into $B_1$. Note that the reduced norm of every element in $B$ is positive 
at all places
in $T$, since $B$ is a matrix algebra over the Hamilton quaternions at those places. Since $L$ embeds 
into $B$, so does the subfield $F$, and in fact, if $\alpha\in F$, then 
$N_{L/K}(\alpha)=N_{F/K}(\alpha)^{[L:F]}>0$ implies
 $N_{F/K}(\alpha)>0$, since $[L:F]$ is odd. We conclude that $F/K$ ramifies at every place in $T$, 
whence $F$ satisfies the
suficient conditions of the Albert-Brauer-Hasse-Noether Theorem 
\cite[Theorem 6.1.1]{linothesis} to embeds into $B_1$.

Now we assume that $L$ satisfies conditions (1)-(4), and let $E/K$ be the quadratic extension in (4).
 It was proved in \cite{smallranks} that $E$, and therefore also $L$,
is selective. Here we give an independent proof for the sake of completeness. Identify $L$ with a subalgebra of $B$.
 Let $\Ha=\oink_L$, and let $\Da$ be a maximal order in $B$ containing $\Ha$. 
It suffices to prove that the representation field $F(\Da/\Ha)$ contains $E$. For this, it suffices to prove that the
relative spinor image $H(\Da|\Ha)$ is contained in $H_E=K^*N_{E/K}(J_E)$.

Note that, by hypothesis, $\oink_E$ is selective in $B_1$, and $N_{B_1}(B_{1\wp}^*)=N_B(B_{\wp}^*)$ 
at every archimedean place, since $B_2$ is a matrix algebra at those places. 
On the other hand, $N_{B_1}(B_{1\wp}^*)$ is contained in $H_E$ at all places, 
since $\oink_E$ is selective for $B_1$. 
It suffices therefore to prove the contention at finite places. Note that $E/K$ is unramified. 
Let $\wp$ be a finite place inert for $E/K$, and let $\mathbb{P}$ be the unique place
of $E$ lying over $\wp$. We need to prove that 
$H_\wp(\Da|\Ha)$ is contained in $\oink_\wp^* K_\wp^{*2}$.  By \cite[Theorem 1]{abelianos}, 
it suffices to prove that every irreducible representation of the residual algebra $\mathbb{H}_\wp$ 
has even dimension. Here,  
$\mathbb{H}_\wp$ is the image of $\Ha$ in the quotient
$\Da_\wp/I$, where $I$ is the unique maximal two-sided ideal of $\Da_\wp$.  
Note that $\mathbb{H}_\wp$ is a quotient
of $\Ha_\wp/\pi\Ha_\wp$ where $\pi$ is a uniformizing parameter of $K_\wp$. 
Since $\Ha$ is an algebra over the ring of integers $\oink_E$,
and therefore $\Ha_\wp$
is an algebra over the ring of integers $\oink_{E_{\mathbb{P}}}$ 
of the unique unramified quadratic extension $E_{\mathbb{P}}$ of $K_\wp$, 
we conclude that  $\mathbb{H}_\wp$ is an algebra over the residue field of $E_{\mathbb{P}}$,
and the result follows.\qed

\end{document}